\documentclass[12pt]{amsart}

\usepackage{mathptmx}
\usepackage{graphicx}
\usepackage{bm}
\usepackage{pstricks}

\usepackage{amsmath}
\usepackage{amsfonts}
\usepackage{amssymb}
\usepackage{amscd}
\usepackage{xypic}
\usepackage{latexsym}
\usepackage{rotating}
\usepackage{times}

\def\ra{\rightarrow}

\def\k#1{\kern -#1pt}

\def\mcc#1{{\tiny\left(\k{5}\begin{array}{cc}#1\end{array}\k5\right)}}

\def\mccc#1{{\tiny\left(\k5\begin{array}{ccc}#1\end{array}\k5\right)}}

\def\mcccc#1{{\tiny\left(\k5\begin{array}{cccc}#1\end{array}\k5\right)}}

\def\mcc#1{{\scriptscriptstyle\left(\k{5}%
\begin{array}{cc}#1\end{array}\k5\right)}}

\def\mccc#1{{\scriptscriptstyle\left(\k{5}%
\begin{array}{ccc}#1\end{array}\k5\right)}}

\def\mcccc#1{{\scriptscriptstyle\left(\k{5}%
\begin{array}{cccc}#1\end{array}\k5\right)}}

\def\A{\mathcal{A}}

\def\C{\mathbb{C}}

\def\N{\mathbb{N}}

\def\HWS{{\cal H}_{WS}}
\def\1H{\mathbf{1}_{\HWS}}
\def\L{\mathbb{L}}
\def\V{\mathbb{V}}

\def\al{\alpha}
\def\be{\beta}

\def\om{\omega}
\def\MQSym{{\bf MQSym}}

\def\2m#1#2{\begin{bmatrix}#1\cr #2\end{bmatrix}}
\def\3m#1#2#3{\begin{bmatrix}#1\cr #2\cr #3\end{bmatrix}}

\def\diag{\mathbf{diag}}
\def\ldiag{\mathbf{ldiag}}
\def\lldiag{\mathbf{lldiag}}
\def\DIAG{\mathbf{DIAG}}
\def\LDIAG{\mathbf{LDIAG}}
\def\MLD{\mathbf{MLDIAG}}
\def\MQS{\mathbf{MQSym}}
\def\SMQ{\mathbf{SMQSym}}
\def\End{\mathrm{End}}

\def\shuff#1#2{\mathbin{
      \hbox{\vbox{
        \hbox{\vrule
              \hskip#2
              \vrule height#1 width 0pt
               }%
        \hrule}%
             \vbox{
        \hbox{\vrule
              \hskip#2
              \vrule height#1 width 0pt
               \vrule }%
        \hrule}%
}}}
\def\shuffl{{\mathchoice{\shuff{7pt}{3.5pt}}%
                        {\shuff{6pt}{3pt}}%
                        {\shuff{4pt}{2pt}}%
                        {\shuff{3pt}{1.5pt}}}}%

\def\ash{\, \underline{\shuffl} \, }

\def\setinterlineskip#1{\baselineskip=0pt
  \lineskip=#1 \lineskiplimit=\maxdimen}
\def\DessinsMatrix#1{\vcenter{\hbox{\makebox[1.3ex]{$\scriptstyle#1$}}}}
\def\GensMatrix#1{\vcenter{\halign{&$\DessinsMatrix{##}$\cr#1}}\egroup}

\def\smallmatrice{%
  \bgroup
  \let\ =\omit
  \let\\=\cr
  \setinterlineskip{8.0pt}
  \GensMatrix}

\def\ss{\smallskip}

\def\bs{\bigskip}

\def\ncp#1#2{#1\langle #2\rangle}

\def\S{\mathfrak S}
\def\A{\mathcal{A}}

\def\C{\mathbb{C}}

\def\N{\mathbb{N}}

\def\ra{\rightarrow}

\def\ss{\smallskip}

\newtheorem{example}{Example}[section]

\newtheorem{theorem}[example]{Theorem}

\newtheorem{corollary}[example]{Corollary}

\newtheorem{proposition}[example]{Proposition}

\newtheorem{lemma}[example]{Lemma}

\def\LDiag{{\bf LDIAG}}

\def\WQSym{{\bf WQSym}}
\def\WSym{{\bf WSym}}
\def\WQ{{\bf WQ}} 
\def\W{{\bf W}} 
\def\FQ{{\bf F}}
\def\F{{\bf F}} 

\def\SMQSym{{\bf SMQSym}}
\def\SMRSym{{\bf SMRSym}}
\def\SMCSym{{\bf SMCSym}}
\def\SMSym{{\bf SMSym}}
\def\SMQ{{\bf SMQ}} 
\def\SMR{{\bf SMR}} 
\def\SMC{{\bf SMC}} 
\def\SM{{\bf SM}} 

\def\MQSym{{\bf MQSym}}
\def\MRSym{{\bf MRSym}}
\def\MCSym{{\bf MCSym}}
\def\MSym{{\bf MSym}}
\def\MQ{{\bf MQ}} 
\def\MR{{\bf MR}} 
\def\M{{\bf M}} 

\def\ssh{\Cup} 
\def\std{{\rm std}}

\def\st{\scriptstyle}
\def\es{\emptyset}

\def\PMuR{{\rm PMuR}}
\def\MuR{{\rm MuR}}

\def\gaudend{\ll}      
\def\droitdend{\gg}    

\def\gf#1#2{\genfrac{}{}{0pt}{}{#1}{#2}}

\title{Hopf algebras of diagrams}

\author[G. H. E. Duchamp, J.-G. Luque, J.-C. Novelli, C. Tollu, F. Toumazet]%
{G. H. E. Duchamp, J.-G. Luque, J.-C. Novelli, C. Tollu, F. Toumazet}

\address[G.H.E. Duchamp, C. Tollu, F. Toumazet]{Institut Galil\'ee, LIPN, CNRS
UMR 7030\\ 99, avenue J.-B. Clement, F-93430 Villetaneuse, France}
\address[J.-G. Luque, J.-C Novelli] {Universite Paris-Est, Institut Gaspard
Monge, \\
5 Boulevard Descartes \\Champs-sur-Marne \\77454 Marne-la-Vall\'ee cedex 2 \\
France}
\email[G. H. E. Duchamp]{ghed@lipn.univ-paris13.fr}
\email[J.-G. Luque]{luque@univ-mlv.fr}
\email[J.-C. Novelli]{novelli@univ-mlv.fr (corresponding author)}
\email[C. Tollu]{ct@lipn.univ-paris13.fr}
\email[F. Toumazet]{ft@lipn.univ-paris13.fr}
\date{\today}



\begin{document}

\begin{abstract}
We investigate several Hopf algebras of diagrams related to Quantum
Field Theory of Partitions and whose product comes from the Hopf algebras
$\WSym$ or $\WQSym$ respectively built on integer set partitions and set
compositions.
Bases of these algebras are indexed either by bipartite graphs (labelled or
unlabbeled) or by packed matrices (with integer or set coefficients).
Realizations on biword are exhibited, and it is shown how these algebras fit
into a commutative diagram. Hopf deformations and dendriform structures
are also considered for some algebras in the picture.
\end{abstract}

\maketitle

\section{Introduction}

\bs The purpose of the present paper is twofold. First, we want to tighten the
links between a body of Hopf algebras related to physics and the realm of
noncommutative symmetric functions, although the latter are no longer
disconnected \cite{Foissy, GOF4, NY}. Second, we aim at providing examples
of combinatorial shifting (a generic way of deforming algebras) and expounding
how the Hopf algebra on packed matrices $\MQS$ could be considered a
construction scheme including its first appearance with integers \cite{DHT} as
a special case.

Our paper is the continuation of \cite{DLPT}, as we go deeper into the
connections between combinatorial Hopf algebras and Feynman diagrams of a
special Field Theory introduced by Bender, Brody and Meister \cite{BBM}. These
Feynman diagrams arise in the expansion of
\begin{equation}
G(z)=\left.\exp
\left\{
  \sum_{n\geq 1}\frac{L_n}{n!}
  \left(z\frac{\partial}{\partial x}\right)^n\right\}
\exp\left\{\sum_{m\geq1}V_m\frac{x^m}{m!}\right\}\right|_{x=0}
\end{equation}
and are bipartite finite graphs with no isolated vertex, and edges weighted
with integers. They are in bijective correspondence with packed matrices of
integers up to a permutation of the columns and a permutation of the rows. The
algorithm constructing the matrix from the associated diagram uses as an
intermediate structure a particular packed matrix whose entries are sets. Such
set matrices appear when one computes the internal product in
$\WSym$~\cite{SR} and in $\WQSym$~\cite{Hiv,NT}, then isomorphic to the
Solomon-Tits algebra. In this context, it becomes natural to investigate
Hopf algebras of (set) packed matrices whose product comes from $\WSym$ or
$\WQSym$.

The paper is organized as follows. In Section \ref{Phys}, the connection
between the Quantum Field Theory of Partitions and a three-parameter
deformation of the Hopf algebra \LDiag\  of labelled diagrams is explained. We
introduce the shifting principle (Subsection \ref{shift}) and give two
illustrations. The first one enables to see $\LDiag$ as the shifted version of
an algebra of unlabelled diagrams.
The second one (Subsection \ref{MQS_k}) explains how to carry over some
constructions from algebras of integer matrices to algebras of set
matrices.\\
In Section \ref{SMQ}, we investigate eight Hopf algebras of
matrices related to labelled or unlabelled diagrams. In particular, we
exhibit realizations on biwords and show how some of these are
bidendriform bialgebras, hence proving those algebras are in particular
self-dual, free and cofree.

{\sc Acknowledgements}.
The first author would like to thank Bodo Lass for an illuminating
seminar talk on the algebraic treatment of bipartite graphs. He is also
greatly indebted to Karol Penson for clearing up the physical origin of the
diagrams.

\section{Hopf algebras coming from physics}\label{Phys}

\subsection{Algebras of diagrams}

Many computations carried out by physicists reduce to the `product formula', a
bilinear coupling between two Taylor expandable functions, introduced by 
C.M. Bender, D.C. Brody, and B.K. Meister in their celebrated 
{\it Quantum field theory of partitions} (henceforth referred to as QFTP)
\cite{BBM}. For an example of such a computation derived from a
partition function linked to the Free Boson Gas model,
see~\cite{GOF9}.

To make the story short, the last expansion of the formula involves
a summation over all diagrams of a certain type \cite{BBM,OPG}, a
labelled version of which is described below. These diagrams
are bipartite graphs with multiple edges.
Bender, Brody and Meister \cite{BBM} introduced QFTP as a toy model
to show that every (combinatorial) sequence of integers can be
represented by Feynman diagrams subject to suited rules.

The case where the expansions of the two functions occurring in the product
formula have constant term 1 is of special interest. The functions can be
then be presented as exponentials which can be regarded as "free" through the
classical Bell polynomials expansion \cite{GOF4} or as coming from the
integration of a Frechet one-parameter group of operators \cite{OPG}. Working
out the formal case, one sees that the coupling results in a summation without
multiplicity of a certain kind of labelled bipartite graphs which are
equivalent, as a data structure, to pairs of unordered partitions of the same
set $\{1,2,...,n\}$.
The sum can be reduced as a sum of topologically inequivalent diagrams (a
monoidal basis of $\DIAG$), at the cost of introducing multiplicities.
Theses graphs, which can be considered as the Feynman diagrams of the QFTP,
generate a Hopf algebra compatible with the product and co-addition on the
multipliers. Interpreting $\DIAG$ as the Hopf homomorphic image of its planar
counterpart, $\LDIAG$, gives access to the noncommutative world and to
deformations: the product is deformed by taking into account, through two
variables, the number of crossings of edges involved in the superposition or
the transposition of two vertices, the coprodut by obtained by interpolating.
This gives the final picture of~\cite{GOF9}.

Labelled diagrams can be identified with their weight functions which
are mappings $\om : \N^+\times\N^+\ra \N$ such that the supporting
subgraph
\begin{equation}
\Gamma_\om=\{(i,j)\in \N^+\times\N^+\ |\ w(i,j)\not=0\}
\end{equation}
has projections \emph{i.e.}, $pr_1(\Gamma_\om)=[1,p];\
pr_2(\Gamma_\om)=[1,q]$ for some $p,q\in \N^+$.

\begin{figure}[ht]
\begin{center}
\includegraphics[scale=0.2]{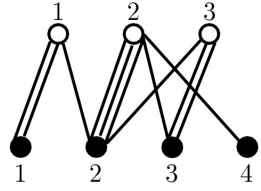}
\caption{\label{labelled} A labelled diagram of shape $3\times 4$.}
\end{center}
\end{figure}

Let $\ldiag$ denote the set of labelled diagrams. With any element $d$ of
$\ldiag$, one can associate the monomial $\L^{\al(d)}\V^{\be(d)}$, called its multiplier, where
$\al(d)$ (resp. $\be(d)$) is the ``white spot type'' (resp. the ``black spot
type'') \emph{i.e.}, the multi-index $(\al_i)_{i\in \N^+}$ (resp.
$(\be_i)_{i\in \N^+}$) such that $\al_i$ (resp.  $\be_i$) is the number of
white spots (resp. black spots) of degree $i$.
For example, the multiplier of the labelled diagram of Figure~\ref{labelled} is
$\L^{(0,0,2,0,1)} \V^{(1,1,1,0,1)}$.

One can endow $\ldiag$ with an algebra structure denoted by $\LDIAG$
where the sum is the formal sum and the product is the shifted concatenation
of diagrams, \emph{i.e.}, consists in juxtaposing the second diagram to the right of the first one
and then adding to the labels of the black spots (resp. of the white spots) of
the second diagram the number of black spots (resp. of white spots) of the
first diagram. Then the application sending a diagram to its multiplier is an
algebra homomorphism.

\begin{figure}[ht]
\begin{center}
\includegraphics[scale=0.4]{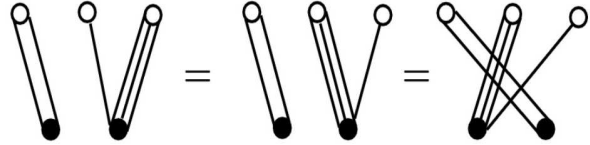}
\caption{\label{nonetiq}Equivalent labelled diagrams.}
\end{center}
\end{figure}

Moreover, the black spots (resp. white spots) of diagram $d$ can be permuted
without changing the monomial $\L^{\al(d)}\V^{\be(d)}$. The classes of
labelled diagrams up to this equivalence relation (permutations of white - or
black - spots among themselves, see Figure~\ref{nonetiq}) are naturally
represented by unlabelled diagrams. The set of unlabelled diagrams will be
henceforth denoted by $\diag$.

The set $\diag$ can also be endowed with an algebra structure, denoted by
$\DIAG$, \emph{e.g.} as the quotient of $\LDIAG$ by the equivalence classes
of labelled diagrams. In $\DIAG$, the product of $d_1$ by $d_2$ is basic
concatenation, {\em i.e.} simply consists in juxtaposing $d_2$ to the
right of $d_1$ \cite{GOF4}.

\subsection{Shifted algebras and applications}
\label{shift}

A three-parameter deformation of the algebra $\LDIAG$ called
$\LDIAG(q_c,q_s,t)$ has been recently constructed, which specializes to both
$\LDIAG$ ($q_c=q_s=t=0$) and $\MQS$ ($q_c=q_s=t=1$).
This construction involves a deformation of the algebra structure, which can
be seen as a particular case of the rather general principle of
\emph{shifting}. This principle will be further exemplified in the sequel of
the present paper.

\begin{lemma} (Shifting lemma.)\\
Let $\A=\oplus_{\al\in M} \A_\al$ be an algebra graded (as a vector space) on
a commutative monoid $(M,+)$.
Let $s\ :\ (M,+)\mapsto (\End_{alg}(\A),\circ)$ be an homomorphism such that
the modified law, given by
\begin{equation}
x *_s y:=x\ s_{\al}(y)  \text{ for all } x\in\A_\al,\ y\in \A
\end{equation}
is $M$-graded. Then, if $\A$ is associative, so is the deformed law $*_s$.
\end{lemma}

Such a procedure, whenever possible, will be called the \emph{shifting of $\A$
by the shift $s$}. We now recall the construction of $\LDIAG(q_c,q_s)$ as
the shifting of another algebra of diagrams.

\ss
Let $M=\N^{(\N^+)}$ be the additive monoid of multidegrees and $M^+=M-\{0\}$
the associated semigroup. A labelled diagram $d$ with $p$ white spots and $q$
black spots can be encoded by a word $W(d)\in (M^+)^*$ as
\begin{equation}
W(d)=\al_1\al_2\cdots \al_q,
\end{equation}
where, for all $i\leq p$, the $i$-th letter of $\al_j$ is the number
of edges joining the black spot $j$ and the white spot $i$
(see Figure~\ref{diag2mot}).
\begin{figure}[h]
\label{diag2mot}
\begin{center}
 \def\JPicScale{0.6}
 \psset{unit=\JPicScale mm}
 \psset{linewidth=0.3,dotsep=1,hatchwidth=0.3,hatchsep=1.5,shadowsize=1,dimen=middle}
 \psset{dotsize=0.7 2.5,dotscale=1 1,fillcolor=black}
 \psset{arrowsize=1 2,arrowlength=1,arrowinset=0.25,tbarsize=0.7 5,bracketlength=0.15,rbracketlength=0.15}
 \begin{pspicture}(-5,65)(233.95,100)
 \rput{0}(58.57,63.5){\psellipse[fillstyle=crosshatch,linewidth=0.5,linestyle=dashed,dash=1 1](0,0)(2.53,2.47)}
 \rput{0}(72.86,63.5){\psellipse*[linewidth=0.75](0,0)(2.53,2.47)}
 \rput{0}(87.14,63.5){\psellipse*[fillstyle=crosshatch,linewidth=0.5,linestyle=dashed,dash=1 1](0,0)(2.53,2.47)}
 \rput{0}(101.43,63.5){\psellipse[fillstyle=crosshatch,linewidth=0.5,linestyle=dashed,dash=1 1](0,0)(2.53,2.47)}
 \rput{0}(60,90){\psellipse[linewidth=0.5](0,0)(2.53,2.47)}
 \rput{0}(80,90){\psellipse[linewidth=0.5](0,0)(2.53,2.47)}
 \rput{0}(99,90){\psellipse[linewidth=0.5](0,0)(2.53,2.47)}
 \psline[linewidth=0.5,linestyle=dashed,dash=1 1](55.75,63.47)(58.57,88)
 \psbezier(58.57,88)(58.57,88)(58.57,88)(58.57,88)
 \psbezier(58.57,88)(58.57,88)(58.57,88)(58.57,88)
 \psline[linewidth=0.5,linestyle=dashed,dash=1 1](57.44,65.68)(57.26,65.68)
 \psline[linewidth=0.5,linestyle=dashed,dash=1 1](60.36,87.97)(57.73,65.71)
 \psline[linewidth=0.75](62.05,87.54)(70.23,63.5)
 \psline(70.23,63.5)(70.23,63.32)
 \psline[linewidth=0.75](70.41,64.7)(78.31,87.82)
 \psline[linewidth=0.75](71.73,65.62)(79.25,87.54)
 \psline[linewidth=0.75](72.86,65.98)(80.28,87.73)
 \psline[linewidth=0.5,linestyle=dashed,dash=1 1](85.45,65.25)(80.94,87.82)
 \psline[linewidth=0.5,linestyle=dashed,dash=1 1](81.97,87.45)(99.55,65.07)
 \psline(74.36,65.34)(74.46,65.34)
 \psline[linewidth=0.75](74.36,65.34)(97.67,87.82)
 \psline[linewidth=0.5,linestyle=dashed,dash=1 1](98.89,87.63)(87.42,65.8)
 \psline[linewidth=0.5,linestyle=dashed,dash=1 1](100.39,87.63)(88.74,65.62)
 \rput{0}(118.57,63.5){\psellipse[fillstyle=crosshatch,linewidth=0.5,linestyle=dashed,dash=1 1](0,0)(2.53,2.47)}
 \rput{0}(132.86,63.5){\psellipse[fillstyle=crosshatch,linewidth=0.5,linestyle=dashed,dash=1 1](0,0)(2.53,2.47)}
 \rput{0}(147.14,63.5){\psellipse*[linewidth=0.75](0,0)(2.53,2.47)}
 \rput{0}(161.43,63.5){\psellipse[fillstyle=crosshatch,linewidth=0.5,linestyle=dashed,dash=1 1](0,0)(2.53,2.47)}
 \rput{0}(120,90){\psellipse[linewidth=0.5,linestyle=dashed,dash=1 1](0,0)(2.53,2.47)}
 \rput{0}(140,90){\psellipse[linewidth=0.5](0,0)(2.53,2.47)}
 \rput{0}(159,90){\psellipse[linewidth=0.5](0,0)(2.53,2.47)}
 \psline[linewidth=0.5,linestyle=dashed,dash=1 1](115.75,63.47)(118.57,88)
 \psbezier[linewidth=0.5,linestyle=dashed,dash=1 1](118.57,88)(118.57,88)(118.57,88)(118.57,88)
 \psbezier[linewidth=0.5,linestyle=dashed,dash=1 1](118.57,88)(118.57,88)(118.57,88)(118.57,88)
 \psline[linewidth=0.5,linestyle=dashed,dash=1 1](117.44,65.68)(117.26,65.68)
 \psline[linewidth=0.5,linestyle=dashed,dash=1 1](120.36,87.97)(117.73,65.71)
 \psline[linewidth=0.5,linestyle=dashed,dash=1 1](122.05,87.54)(130.23,63.5)
 \psline[linewidth=0.5,linestyle=dashed,dash=1 1](130.23,63.5)(130.23,63.32)
 \psline[linewidth=0.5,linestyle=dashed,dash=1 1](130.41,64.7)(138.31,87.82)
 \psline[linewidth=0.5,linestyle=dashed,dash=1 1](131.73,65.62)(139.25,87.54)
 \psline[linewidth=0.5,linestyle=dashed,dash=1 1](132.86,65.98)(140.28,87.73)
 \psline[linewidth=0.75](145.45,65.25)(140.94,87.82)
 \psline[linewidth=0.5,linestyle=dashed,dash=1 1](141.97,87.45)(159.55,65.07)
 \psline[linewidth=0.5,linestyle=dashed,dash=1 1](134.36,65.34)(134.46,65.34)
 \psline[linewidth=0.5,linestyle=dashed,dash=1 1](134.36,65.34)(157.67,87.82)
 \psline[linewidth=0.75](158.89,87.63)(147.42,65.8)
 \psline[linewidth=0.75](160.39,87.63)(148.74,65.62)
 \rput{0}(178.57,63.5){\psellipse[fillstyle=crosshatch,linewidth=0.5,linestyle=dashed,dash=1 1](0,0)(2.53,2.47)}
 \rput{0}(192.86,63.5){\psellipse[fillstyle=crosshatch,linewidth=0.5,linestyle=dashed,dash=1 1](0,0)(2.53,2.47)}
 \rput{0}(207.14,63.5){\psellipse[fillstyle=crosshatch,linewidth=0.5,linestyle=dashed,dash=1 1](0,0)(2.53,2.47)}
 \rput{0}(221.43,63.5){\psellipse*[linewidth=0.75](0,0)(2.53,2.47)}
 \rput{0}(180,90){\psellipse[linewidth=0.5,linestyle=dashed,dash=1 1](0,0)(2.53,2.47)}
 \rput{0}(200,90){\psellipse[linewidth=0.5](0,0)(2.53,2.47)}
 \rput{0}(219,90){\psellipse[linewidth=0.5,linestyle=dashed,dash=1 1](0,0)(2.53,2.47)}
 \psline[linewidth=0.5,linestyle=dashed,dash=1 1](175.75,63.47)(178.57,88)
 \psbezier[linewidth=0.5,linestyle=dashed,dash=1 1](178.57,88)(178.57,88)(178.57,88)(178.57,88)
 \psbezier[linewidth=0.5,linestyle=dashed,dash=1 1](178.57,88)(178.57,88)(178.57,88)(178.57,88)
 \psline[linewidth=0.5,linestyle=dashed,dash=1 1](177.44,65.68)(177.26,65.68)
 \psline[linewidth=0.5,linestyle=dashed,dash=1 1](180.36,87.97)(177.73,65.71)
 \psline[linewidth=0.5,linestyle=dashed,dash=1 1](182.05,87.54)(190.23,63.5)
 \psline[linewidth=0.5,linestyle=dashed,dash=1 1](190.23,63.5)(190.23,63.32)
 \psline[linewidth=0.5,linestyle=dashed,dash=1 1](190.41,64.7)(198.31,87.82)
 \psline[linewidth=0.5,linestyle=dashed,dash=1 1](191.73,65.62)(199.25,87.54)
 \psline[linewidth=0.5,linestyle=dashed,dash=1 1](192.86,65.98)(200.28,87.73)
 \psline[linewidth=0.5,linestyle=dashed,dash=1 1](205.45,65.25)(200.94,87.82)
 \psline[linewidth=0.75](201.97,87.45)(219.55,65.07)
 \psline[linewidth=0.5,linestyle=dashed,dash=1 1](194.36,65.34)(194.46,65.34)
 \psline[linewidth=0.5,linestyle=dashed,dash=1
1](194.36,65.34)(217.67,87.82)
 \psline[linewidth=0.5,linestyle=dashed,dash=1
1](218.89,87.63)(207.42,65.8)
 \psline[linewidth=0.5,linestyle=dashed,dash=1
1](220.39,87.63)(208.74,65.62)
 \rput{0}(-1.43,63.5){\psellipse*[linewidth=0.75](0,0)(2.53,2.47)}
 \rput{0}(12.86,63.5){\psellipse[fillstyle=crosshatch,linewidth=0.5,linestyle=dashed,dash=1
1](0,0)(2.53,2.47)}
 \rput{0}(27.14,63.5){\psellipse[fillstyle=crosshatch,linewidth=0.5,linestyle=dashed,dash=1
1](0,0)(2.53,2.47)}
 \rput{0}(41.43,63.5){\psellipse[fillstyle=crosshatch,linewidth=0.5,linestyle=dashed,dash=1
1](0,0)(2.53,2.47)}
  \rput{0}(0,90){\psellipse[linewidth=0.5](0,0)(2.53,2.47)}
 \rput{0}(20,90){\psellipse[linewidth=0.5,linestyle=dashed,dash=1 1](0,0)(2.53,2.47)}
 \rput{0}(39,90){\psellipse[linewidth=0.5,linestyle=dashed,dash=1 1](0,0)(2.53,2.47)}
 \psline[linewidth=0.75](-4.25,63.47)(-1.43,88)
 \psbezier(-1.43,88)(-1.43,88)(-1.43,88)(-1.43,88)
 \psbezier(-1.43,88)(-1.43,88)(-1.43,88)(-1.43,88)
 \psline(-2.56,65.68)(-2.74,65.68)
 \psline[linewidth=0.75](0.36,87.97)(-2.27,65.71)
 \psline[linewidth=0.5,linestyle=dashed,dash=1
1](2.05,87.54)(10.23,63.5)
 \psline[linewidth=0.5,linestyle=dashed,dash=1
1](10.23,63.5)(10.23,63.32)
 \psline[linewidth=0.5,linestyle=dashed,dash=1
1](10.41,64.7)(18.31,87.82)
 \psline[linewidth=0.5,linestyle=dashed,dash=1
1](11.73,65.62)(19.25,87.54)
 \psline[linewidth=0.5,linestyle=dashed,dash=1
1](12.86,65.98)(20.28,87.73)
 \psline[linewidth=0.5,linestyle=dashed,dash=1
1](25.45,65.25)(20.94,87.82)
 \psline[linewidth=0.5,linestyle=dashed,dash=1
1](21.97,87.45)(39.55,65.07)
 \psline[linewidth=0.5,linestyle=dashed,dash=1
1](14.36,65.34)(14.46,65.34)
 \psline[linewidth=0.5,linestyle=dashed,dash=1
1](14.36,65.34)(37.67,87.82)
 \psline[linewidth=0.5,linestyle=dashed,dash=1
1](38.89,87.63)(27.42,65.8)
 \psline[linewidth=0.5,linestyle=dashed,dash=1
1](40.39,87.63)(28.74,65.62)
 \end{pspicture}

¨
\end{center}
\begin{flushleft}
The edges adjacent to the blackspots correspond successively to the
multidegrees $(2)$ (two edges to the first white spot), $(1,3,1)$
(one edge to the first and third white spots and three edges to the second
one), $(0,1,2)$ and $(0,1)$.
Thus the code is $W(d)=(2)(1,3,1)(0,1,2)(0,1)$.\end{flushleft}
\caption{Coding a diagram with a word of multidegrees. }
\end{figure}

\noindent
One can skew the product in $\ncp{k}{M^+}$, counting crossings and
superpositions as for $\LDIAG$.

\begin{proposition}
Let $k$ be a ring, $q_c,q_s\in k$ and consider the deformed graded law
defined on $\ncp{k}{M^+}$ by

\begin{equation}
\left\{
\begin{array}{rcl}
1_{(M^+)^*}*w & = & w*1_{(M^+)^*}=w,\\[5pt]
\al u * \be v & = & \al(u*\be v)+ q_c^{|\al u||\be|} \be (\al u *v) +
q_s^{|\al||\be|}q_c^{|u||\be|} (\al+\be)(u*v),
\end{array}
\right.
\end{equation}
where $\alpha,\beta\in M^+$, $u,v\in(M^+)^*$, the weight
$|\alpha|$ of multidegree $\alpha$ is just the sum of its
coordinates and the weight $|u|$ of the word
$u=\alpha_1\cdots\alpha_t$ is $|u|=\sum_{i=1}^t|\alpha_i|$.

This product is associative.
\end{proposition}

The algebra $(k\langle M^+\rangle,+,*)$ is denoted by $\MLD (q_c,q_s)$. In the
shifted version, the product amounts to performing all superpositions of black
spots and/or crossings of edges, weighting them with the corresponding value
$q_s$ or $q_c$, powered by the number of crossings of edges.

This construction is reminiscent, up to the deformations, of
Hoffman's~\cite{Ho} and its variants~\cite{Ca,FG}, and also of an older
one, the infiltration product in computer science~\cite{Oc,DFLL}.

\ss The shift going from $\MLD(q_c,q_s)$ to itself is the following.
Let $\al_1\al_2\cdots \al_p\in (M^+)^*$ and $n\in \N$. One sets
\begin{equation}
s_n(\al_1\al_2\cdots \al_p):=(0^n\al_1)(0^n\al_2)\cdots (0^n\al_p)
\end{equation}
where $0^n\al$ is the insertion of $n$ zeroes on the left of $\al$.
Note that $n\mapsto s_n$ is a homomorphism of monoids $(\N,+)\ra
\End_{alg}\big(\MLD (q_c,q_s)\big)$.

\subsection{Another application of the shifting principle}
\label{MQS_k}
The Hopf operations of $\MQS$ as described in~\cite{DHT} do not depend on the
fact that the entries of the matrices are integers.

For a pointed set $(X,x_0)$ (\emph{i.e.}, $x_0\in X$), let us denote by
$\MQS_k(X,x_0)$ the $k$-vector space spanned by rectangular matrices with
entries in $X$ with no line or column filled with $x_0$ (which plays now the
r\^ole of zero) and the product, coproduct, unit and counit as in~\cite{DHT}.
It is clear that $\MQS_k(X,x_0)$ is a Hopf algebra and that the correspondence
$(X,x_0)\mapsto \MQS_k(X,x_0)$ is a functor from the category of pointed sets
(endowed with the strict arrows, that is, the mappings $\phi\ :\ (X,x_0)\ra
(Y,y_0)$ such that $\phi(x_0)=y_0$ and $\phi(X-\{x_0\})\subset Y-\{y_0\}$) to
the category of $k$-Hopf algebras.
In the particular case when $X=2^{(\N^+)}$, that is, finite subsets of $\N^+$,
and $x_0=\emptyset$, one can define a shift by a translation of the elements.
More precisely, for $F\in 2^{(\N^+)}$ and $M$ a $p\times q$ matrix
with coefficients in $2^{(\N^+)}$, one sets
\begin{equation}
s_n(F)=\{x+n\}_{x\in F}\ ,\ s_n(M)
   =\big(s_n(M[i,j])\big)_{\genfrac{}{}{0pt}{}{1\leq i\leq p}{1\leq j\leq q}}
\,.
\end{equation}

One can check that the $s_n$ define a shift on $\MQS_k(2^{(N^+)},\emptyset)$
for the grading given by
\begin{equation}
\MQS_k(2^{(N^+)},\emptyset)_n=
\textrm{span of the matrices whose maximum entry is }n.
\end{equation}

For example, the vector space generated by the packed matrices whose
entries partition the set $\{1,2,\cdots n\}$ is closed by the product. This
is the algebra $\SMQSym$. We shall see that this corresponds to labelling
the edges of a labelled diagram with numbers from 1 to $k$.

\section{Packed matrices and related Hopf algebras}
\label{SMQ}
\subsection{The combinatorial objects}

In the sequel, we represent different kinds of diagrams using matrices to
emphasize the parallel between this construction and the construction of
$\MQSym$ (\cite{DHT}).

\subsubsection{Set packed matrices}

Since the computations are the same in many cases, let us begin with the most
general case and explain how one recovers the other cases by algebraic means.
Let us consider the set $\lldiag$ of bipartite graphs with white and black
vertices, and edges, all three labelled by initial intervals $[1,p]$
of $\N^*$.
The diagrams $\ldiag$ are obtained by erasing the labels of the edges of such
an element.


The set $\lldiag$ is in direct bijection with \emph{set packed matrices}, that
are matrices containing disjoint subsets of $[1,n]$ for some $n\in\N$ with
no line or column filled with empty sets, and such that the union of all
subsets is $[1,n]$ itself. The bijection consists in putting $k$ in the
$(i,j)$ entry of the matrix if the edge labelled $k$ connects the white dot
labelled $i$ with the black dot labelled $j$.
Figure~\ref{mat-sc} shows an example of such a matrix.

\begin{figure}[ht]
\begin{equation*}
\left(
 \begin{array}{ccc}
  \{3\}    &\{6\} &\{2\}\\
  \emptyset&\{1,5\}   &\{4\}
  \end{array}\right)
\end{equation*}
\caption{\label{mat-sc}A set packed matrix.}
\end{figure}

Note that set packed matrices are in bijection with pairs of
\emph{set compositions}, or, ordered set partitions of $[1,n]$: given a set
packed matrix, compute the ordered sequence of the union of the elements in
the same row (resp. column).
For example, the set packed matrix of Figure~\ref{mat-sc} gives rise to the
set compositions $[\{2,3,6\},\{1,4,5\}]$ and $[\{3\},\{1,5,6\},\{2,4\}]$.
Given two set compositions $\Pi$ and $\Pi'$, define
$M_{ij}:=\Pi_i\cap\Pi'_j$.

So the generating series counting set packed matrices by their maximum entry
$n$ is given by the square of the ordered Bell numbers, that is
sequence A122725 in~\cite{Slo}.

\subsubsection{Integer packed matrices}

As already said, if one forgets the labels of the edges of an element of
$\lldiag$, one recovers an element of $\ldiag$. Its matrix representation
is an \emph{integer packed matrix}, that is, a matrix with no line or column
filled with zeros. The encoding is simple: $m_{ij}$ is equal to the number of
edges between the white spot labelled $i$ and the black spot labelled $j$.
Note that, from the matrix point of view, it consists in replacing the subsets
by their cardinality.
Figure~\ref{mat-int} shows an example of such a matrix.

\begin{figure}[ht]
\begin{equation*}
\left(
 \begin{array}{ccc}
 1&1&1\\
 0&2&1\\
 \end{array}\right)
\end{equation*}
\caption{\label{mat-int}An integer packed matrix.}
\end{figure}

The generating series counting integer packed matrices by the sum of their
entries is given by sequence A120733 of~\cite{Slo}.

\subsubsection{Other packed matrices}

In the sequel, we shall also consider diagrams where one forgets about the
labels of the white spots, or about the labels of the black spots, or about
all labels.
Those three classes of diagrams are respectively in bijection with matrices up
to a permutation of the rows, a permutation of the columns, and simultaneous
permutations of both.

\subsection{Word quasi-symmetric and symmetric functions}

Let us recall briefly the definition of two combinatorial Hopf algebras that
will be useful in the sequel.

\subsubsection{The Hopf algebra $\WQSym$}

We use the notations of \cite{NT}. The word quasi-symmetric
functions are the noncommutative polynomial invariants of Hivert's
quasi-symmetrizing action \cite{Hiv}
\begin{equation}
\WQSym(A):=\C\langle A \rangle^{{\mathfrak S}(A)_{QS}}.
\end{equation}

When $A$ is an infinite alphabet, $\WQSym(A)$ is a graded Hopf algebra whose
basis is indexed by set compositions, or, equivalently, packed words. Recall
that packed words are words $w$ on the alphabet $[1,k]$ such that if $i\not=1$
appears in $w$, then $i-1$ also appears in $w$. The bijection between both
sets is that $w_i=j$ iff $i$ is in the $j$-th part of the set composition (see
Figure~\ref{sc-pack} for an example).

\begin{figure}[ht]
\begin{center}
$cbbadacb \qquad\longleftrightarrow\qquad [\{4,6\},\{2,3,8\},\{1,7\},\{5\}]$
\end{center}
\caption{\label{sc-pack}A packed word and its corresponding set composition.}
\end{figure}

By definition, $\WQSym$ is generated by the polynomials
\begin{equation}
\WQ_u:=\sum_{{\rm pack}(w)=u}w,
\end{equation}
where $u={\rm pack}(w)$ is the packed word having the same comparison
relations between all elements as $w$.

The product in $\WQSym$ is given by
\begin{equation}
\WQ_u\WQ_v=\sum_{w\in u\star_Wv}\WQ_w,
\end{equation}
where the convolution $u\star_Wv$ of two packed words is defined by
\begin{equation}
u\star_Wv :=
\sum_{\gf{w_1,w_2; w=w_1w_2\ {\rm packed}}%
         {{\rm pack}(w_1)=u,\ {\rm pack}(w_2)=v}}
      w.
\end{equation}
The coproduct is given by
\begin{equation}
\Delta \WQ_w(A)= \sum_{u,v ; w\in u\ssh_W\,v} \WQ_u\otimes \WQ_v
\end{equation}
where $u\ssh_W\, v$ denotes the \emph{packed shifted shuffle} that
is the shuffle of $u$ and $v'=v[max(u)]$, that is, $v'_i=v_i+max(u)$.

The dual algebra $\WQSym^*$ of $\WQSym$ is a subalgebra of the Parking
quasi-symmetric functions $\bf PQSym$~\cite{NT-bi}.
This algebra has a multiplicative basis denoted by $\FQ^w$, where the product
is the shifted concatenation, that is $u. v[max(u)]$.

\subsubsection{$\WSym$}

The algebra of word symmetric functions $\WSym$, first defined by Rosas and
Sagan in~\cite{SR}, where it is called the algebra of symmetric functions in
noncommuting variables, is the Hopf subalgebra of $\WQSym$ generated by
\begin{equation}
\W_{\pi}:=\sum_{\mbox{\scriptsize sp}(u)=\pi}\WQ_{u},
\end{equation}
where $\mbox{sp}(u)$ is the (unordered) set partition obtained by forgetting the
order of the parts of its corresponding set composition.

Its dual $\WSym^*$ is the quotient of $\WQSym^*$
\begin{equation}
\WSym^*=\WQSym^*/_J
\end{equation}
where $J$ is the ideal generated by the polynomials $\FQ^u-\FQ^v$
with $u$ and $v$ corresponding to the same set partition.
We denote by $\F^{\mbox{\scriptsize sp}(u)}$ the image of $\FQ^u$ by the canonical surjection.

\section{Hopf algebras of set packed matrices}

\subsection{Set matrix quasi-symmetric functions}

The construction of the Hopf algebra $\SMQSym$ over set packed matrices is a
direct adaptation of the construction of $\MQSym$ (\cite{Hiv,DHT}).
Consider the linear subspace spanned by the elements $\SMQ_{M}$, where $M$
runs over the set of set packed matrices. We denote by $h(M)$ the number of
rows of $M$. Then define
\begin{equation}
\SMQ_{P} \SMQ_{Q} := \sum_{R\in\ash (P,Q)} \SMQ_R
\end{equation}

where the {\em augmented shuffle} of $P$ and $Q$,
$\ash(P,Q)$ is defined as follows: let $Q'$ be obtained from $Q$ by adding the
greatest number inside $P$ to all elements inside $Q$. Let $r$ be an
integer between $\max(p,q)$ and $p+q$, where $p=h(P)$ and $q=h(Q)$.
Insert rows of zeros in the matrices $P$ and $Q'$ so as to form
matrices $\tilde P$ and $\tilde Q'$ of height $r$. Let $R$ be
the matrix obtained by gluing $\tilde Q'$ to the right of $\tilde P$.
The set $\ash (P,Q)$ is formed by all matrices with no row of 0's obtained
this way.

For example,
\begin{equation}
\begin{array}{l}
  \SMQ_{\mcc{\st\{3,4\}&\k8\st\{1\}\\\st\{2\}&\k8\st\emptyset}}
  \SMQ_{\mcc{\st\{2,3,4\}&\k8\st\{1\}}} =

  \SMQ_{\mcccc{
             \st\{3,4\}&\k8\st\{1\}&\k8\st\es&\k8\st\es\\
             \st\{2\}&\k8\st\es&\k8\st\es&\k8\st\es\\
             \st\es&\k8\st\es&\k8\st\{6,7,8\}&\k8\st\{5\}}}
 +\SMQ_{\mcccc{
             \st\{3,4\}&\k8\st\{1\}&\k8\st\es&\k8\st\es\\
             \st\{2\}&\k8\st\es&\k8\st\{6,7,8\}&\k8\st\{5\}}}
 \\[3mm]
\qquad\
  +\SMQ_{\mcccc{
             \st\{3,4\}&\k8\st\{1\}&\k8\st\es&\k8\st\es\\
             \st\es&\k8\st\es&\k8\st\{6,7,8\}&\k8\st\{5\}\\
             \st\{2\}&\k8\st\es&\k8\st\es&\k8\st\es}}
  +\SMQ_{\mcccc{
             \st\{3,4\}&\k8\st\{1\}&\k8\st\{6,7,8\}&\k8\st\{5\}\\
            \st\{2\}&\k8\st\es&\k8\st\es&\k8\st\es}}
  +\SMQ_{\mcccc{
             \st\es&\k8\st\es&\k8\st\{6,7,8\}&\k8\st\{5\}\\
             \st\{3,4\}&\k8\st\{1\}&\k8\st\es&\k8\st\es\\
             \st\{2\}&\k8\st\es&\k8\st\es&\k8\st\es}}.
\end{array}
\end{equation}

The coproduct $\Delta \SMQ_{M}$ is defined by
\begin{equation}
\Delta\SMQ_A=\sum_{A=\binom{A_1}{A_2}}
                 \SMQ_{\std(A_1)}\otimes \SMQ_{\std(A_2)},
\end{equation}
where $\std(A)$ denotes the standardized of the matrix $A$, that is the matrix
obtained by the substitution $a_i\mapsto i$, where $a_1<\cdots<a_n$ are the integers
appearing in $A$. For example,

\begin{equation}
\begin{split}
\Delta\SMQ_{\mcc{ \st\{2,4\}&\k8\st\{1\}\\
                                    \st\{6\}&\k8\st\{3,5\}}}
=& \SMQ_{\mcc{ \st\{2,4\}&\k8\st\{1\}\\
                                 \st\{6\}&\k8\st\{3,5\}}}\otimes
1+ \SMQ_{\mcc{ \st\{2,3\}&\k8\st\{1\}}}\otimes
   \SMQ_{\mcc{\st\{3\}&\k8\st\{1,2\}}}
\\
&+1\otimes
\SMQ_{\mcc{\st\{2,4\}&\k8\st\{1\}\\
                                    \st\{6\}&\k8\st\{3,5\}}}
\,.
\end{split}
\end{equation}

Rather than checking the compatibility between the product and the coproduct,
one can look for a realization of $\SMQSym$, here in terms of noncommutative
bi-words, that will later give useful guidelines to select or understand
homomorphisms between the different algebras.

\subsection{Realization of $\SMQSym$}

A noncommutative bi-word is a word over an alphabet of bi-letters
$\left\langle \gf{a_i}{{\bf b}_j}\right\rangle$ with $i,j\in\N$ where $a_i$
and ${\bf b}_j$ are letters of two distinct ordered alphabets $A$ and $\bf B$.
We will denote by $\C\left\langle \gf{A}{{\bf B}}\right\rangle$ the algebra of
the (potentially infinite) polynomials over the bi-letters
$\left\langle \gf{a}{{\bf b}}\right\rangle$ for the product $\star$ defined by
\begin{equation}
\left\langle \gf{u_1}{{\bf v}_1}\right\rangle
\star\left\langle \gf{u_2}{{\bf v}_2}\right\rangle=
\left\langle \gf{u_1.u_2}{{\bf v}_1.{\bf v}_2[\max({\bf v}_1)]}\right\rangle,
\end{equation}
where ${\bf v}_1.{\bf v}_2[k]$ denotes the concatenation of ${\bf v}_1$ with
the word ${\bf v}_2$ whose letters are shifted by $k$ and $\max({\bf v}_1)$
denotes the maximum letter of ${\bf v}_1$.

In the sequel, we shall forget the letters $a$ and $\bf b$ when there is no
ambiguity about the alphabets $A$ and $\bf B$, so that, for example, 
\begin{equation}
\left\langle \gf{1\,4\,2}{2\,3\,6}\right\rangle
\star \left\langle \gf{2\,4}{3\,1}\right\rangle
=\left\langle \gf{1\,4\,2\,2\,4}{2\,3\,6\,9\,7}\right\rangle.
\end{equation}

\begin{lemma}
The product $\star$ is associative.
\end{lemma}

\begin{proof}
Straightforward from its definition.
\end{proof}

With each set packed matrix, one associates the bi-word whose $i$th
bi-letter is the coordinate in which the letter $i$ appears in the matrix.
For example,

\begin{equation}
\mbox{bi-word}
\left(\begin{array}{cccc}
\{2,7\}&\es&\es&\{3,5\}\\
\{8\}&\es&\es&\es\\
\es&\{1\}&\{4,6\}&\es\\
\end{array}\right)
= \left\langle \gf{3\,1\,1\,3\,1\,3\,1\,2}{2\,1\,4\,3\,4\,3\,1\,1}
  \right\rangle.
\end{equation}

A bi-word is said bi-packed if its two words are packed.
The bi-packed of a bi-word is the bi-word obtained by packing its
two words.
The set packed matrices are obviously in bijection with the
bi-packed bi-words.

\begin{theorem}
\label{Endgr}
Let $\left\langle \gf{u}{{\bf v}}\right\rangle$
be a bi-packed bi-word. Then
\begin{itemize}
\item The algebra $\SMQSym$ can be realized on bi-words by
\begin{equation}
\SMQ_{\left\langle \gf{u}{{\bf v}}\right\rangle}:=
\sum_{{\rm bipacked }\left\langle \gf{u'}{{\bf v}'}\right\rangle
      =\left\langle \gf{u}{{\bf v}}\right\rangle
     } \left\langle \gf{u'}{{\bf v'}}\right\rangle.
\end{equation}
\item $\SMQSym$ is a Hopf algebra.

\item $\SMQSym$ is isomorphic as a Hopf algebra to the graded endomorphisms of
$\WQSym$:
\begin{equation}
{\rm End}_{gr}\WQSym=\bigoplus_n \WQSym_n\otimes \WQSym^*_n
\end{equation}
through the Hopf homomorphism
\begin{equation}
\phi\left(\SMQ_{\left\langle \gf{u}{\bf v}\right\rangle}\right)
= \WQ_{u}\otimes \FQ^{\bf v}.
\end{equation}
\end{itemize}
\end{theorem}

\begin{proof}
Since the map sending each set packed matrix to a bi-packed bi-word is a
bijection, the first part of the theorem amounts to checking the compatibility
of the product,
\begin{equation}
{\bf SMQ}_{{\rm bi-word}(P)}{\bf SMQ}_{{\rm
bi-word}(Q)}=\sum_{R\in\ash (P,Q)}{\bf SMQ}_{{\rm bi-word}(R)},
\end{equation}
which is straightforward from the definition.
%

${\bf SMQSym}$ being a connected graded algebra, it suffices to
show that the coproduct $\Delta$ is a homomorphism of algebras.
If one uses the representation of the basis elements by pairs of set
compositions, the coproduct reads
\begin{equation}
\Delta\SMQ_{\Pi_1,\Pi_2}=\sum_{\Pi_1=[\Pi'_1,\Pi_1'']}
\SMQ_{\std(\Pi'_1),\std(\Pi_2|_{\Pi'_1})}\otimes
\SMQ_{\std(\Pi''_1),\std(\Pi_2|_{\Pi''_1})},
\end{equation}
where $\Pi_1|_{\Pi_2}$ is the list of sets $[(\Pi_1)_i\cap
\bigcup_j(\Pi_2)_j]_i$ from which one erases the empty sets. The
proof then amounts to mimicking the proof that $\WQSym$ is a Hopf algebra
(see \cite{Hiv}).

One endows ${\rm End}_{gr}\WQSym$ with the coproduct $\Delta$ defined by
\begin{equation}
\Delta {\bf WQ}_{\Pi_1}\otimes {\bf F}^{\Pi_2}=\sum_{\Pi_1=[\Pi'_1,\Pi''_1]}
({\bf WQ}_{\Pi'_1}\otimes {\bf F}^{\Pi_2|_{\Pi'_1}})\otimes ({\bf
WQ}_{\Pi''_1}\otimes {\bf F}^{\Pi_2|_{\Pi''_1}}).
\end{equation}
One then easily checks that $\phi$ is a surjective Hopf homomorphism and since
the two spaces have same series of dimensions, we get the result.
\end{proof}

For example,
\begin{equation}
\SMQ_{\mcccc{
             \st\{3,4\}&\k8\st\{1\}&\k8\st\es&\k8\st\es\\
             \st\es&\k8\st\es&\k8\st\{6\}&\k8\st\es\\
             \st\{2\}&\k8\st\es&\k8\st\es&\k8\st\{5,7\}}}
= \SMQ_{\left\langle \gf{1\,3\,1\,1\,3\,2\,3}{2\,1\,1\,1\,4\,3\,4}
        \right\rangle}
= \sum_{ \gf{j_1<j_2<j_3}{k_1<k_2<k_3<k_4}}
        \left\langle
        \gf{j_1\,j_3\,j_1\,j_1\,j_3\,j_2\,j_3}%
           {k_2\,k_1\,k_1\,k_1\,k_4\,k_3\,k_4} \right\rangle.
\end{equation}
Note that, from the point of view of the realization, the coproduct of
$\SMQSym$ is given by the usual trick of noncommutative symmetric
functions, considering an alphabet $A$ of bi-letters ordered lexicographically
as an ordered sum of two mutually commuting alphabets $A'\hat+A''$ of
bi-letters such that if $(x,y)$ is in $A$ then so is any bi-letter of the
form $(x,z)$. Then the coproduct is a homomorphism for the product.

\subsection{Set matrix half-symmetric functions}

\subsubsection{The Hopf algebra $\SMRSym$}

Let $\SMRSym$ be the subalgebra of $\SMQSym$ generated by the polynomials
$\SMR_{\pi_1,\Pi_2}$ indexed by a set partition $\pi_1$ and a set composition
$\Pi_2$ and defined by

\begin{equation}
\SMR_{(\pi_1,\Pi_2)}
:=\sum_{{\rm sp}(\Pi_1)=\pi_1} \SMQ_{\Pi_1,\Pi_2}.
\end{equation}

For example,
\begin{equation}
\begin{array}{rcl}
\SMR_{\{\{14\},\{2\},\{3\}\},[\{134\},\{2\}]}&=&
  \SMQ_{\mcc{\st\{1,4\}&\k8\st\es\\\st\es&\k8\st\{2\}\\\st\{3\}&\k8\st\es}}
+ \SMQ_{\mcc{\st\{1,4\}&\k8\st\es\\\st\{3\}&\k8\st\es\\\st\es&\k8\st\{2\}}}
+ \SMQ_{\mcc{\st\es&\k8\st\{2\}\\\st\{1,4\}&\k8\st\es\\\st\{3\}&\k8\st\es}}\\
&&
+ \SMQ_{\mcc{\st\es&\k8\st\{2\}\\\st\{3\}&\k8\st\es\\\st\{1,4\}&\k8\st\es}}
+ \SMQ_{\mcc{\st\{3\}&\k8\st\es\\\st\es&\k8\st\{2\}\\\st\{1,4\}&\k8\st\es}}
+ \SMQ_{\mcc{\st\{3\}&\k8\st\es\\\st\{1,4\}&\k8\st\es\\\st\es&\k8\st\{2\}}}.
 \end{array}
\end{equation}

Note that a pair constituted by a set partition and a set composition is
equivalent to a set packed matrix up to a permutation of its rows. Hence, the
realization on bi-words follows: for example,

\begin{equation}
{\rm \SMR}_{\{\{14\},\{2\},\{3\}\},[\{134\},\{2\}]}
=\sum_{\gf{j_1, j_2, j_3\text{ distinct}}{k_1<k_2}}
 \left\langle \gf{j_1\,j_2\,j_3\,j_1}{k_1\,k_2\,k_1\,k_1}
 \right\rangle\,.
\end{equation}

\begin{proposition}
~

\begin{itemize}
\item $\SMRSym$ is isomorphic to
 $\oplus_n\left( \WSym_n\otimes \WQSym_n^*\right)$.
\item $\SMRSym$ is a co-commutative Hopf subalgebra of $\SMQSym$.
\end{itemize}
\end{proposition}

\begin{proof}
The first part of the proposition is a direct consequence of the following
sequence of equalities:
\begin{equation}
\begin{array}{rcl}
\phi(\SMR_{\pi_1,\Pi_2})
&=&\sum_{{\rm sp}(\Pi_1)=\pi_1}\phi({\rm \SMQ}_{\Pi_1,\Pi_2})\\[.2cm]
&=&\sum_{{\rm sp}(\Pi_1)=\pi_1} \WQ_{\Pi_1}\otimes \FQ^{\Pi_2}\\[.2cm]
&=& \W_{\pi_1}\otimes \FQ^{\Pi_2}.\end{array}
\end{equation}

From its definition, $\SMRSym$ is stable for the product and $\Delta$ maps
$\SMRSym$ to $\SMRSym\otimes\SMRSym$. It follows that $\SMRSym$ is a  Hopf
subalgebra of $\SMQSym$. One checks easily the co-commutativity by
restricting $\Delta$ to $\SMRSym$.
\end{proof}

\subsubsection{The Hopf algebra $\SMCSym$}

Forgetting about the order of the columns instead of the rows leads to another
Hopf algebra, $\SMCSym$, a basis of which is indexed by pairs
$(\Pi_1,\pi_2)$ where $\Pi_1$ is a set composition and $\pi_2$ is a
set partition.
It is naturally the quotient (and not a subalgebra) of $\SMQSym$ by the ideal
generated by the polynomials
\begin{equation}
{\rm \SMQ}_{\Pi_1,\Pi_2}-{\rm \SMQ}_{\Pi_1,\Pi'_2}
\end{equation}
where ${\rm sp}(\Pi_2)={\rm sp}(\Pi'_2)$.
Note that this quotient can be brought down to the bi-words.
We denote by $\alpha$ the canonical surjection:
\begin{equation}
\alpha(\SMQ_{\Pi_1,\Pi_2})=:{\rm \SMC}_{\Pi_1,{\rm sp}(\Pi_2)}
\end{equation}

\begin{proposition}
\label{PSMCSym}
~

\begin{itemize}
\item $\SMCSym$ is isomorphic to $\oplus_n \WQSym_n\otimes \WSym_n^*$,
\item $\SMCSym$ is a Hopf algebra.
\end{itemize}
\end{proposition}

\begin{proof}
The first property follows from the fact that the following diagram is
commutative:
\begin{equation}
\begin{array}{ccc}
\SMQSym&\displaystyle\mathop{\rightarrow}^{\phi}&
\oplus_n \WQSym_n\otimes \WQSym_n^*\\
\downarrow\alpha&&Id\otimes\alpha'\downarrow\\
\SMCSym&\displaystyle\mathop{\rightarrow}^{\phi^C}&
\oplus_n \WQSym_n\otimes \WSym_n^*,
\end{array}
\end{equation}
where $\alpha'$ denotes the canonical surjection
$\alpha':\WQSym_n^*\rightarrow \WSym_n^*$, and $\phi^C$ is the map
sending $\SMC_{\Pi_1,\pi_2}$ to $\W_{\Pi_1}\otimes F^{\pi_2}$.
Indeed, the image by $\phi$ of the ideal generated by the
polynomials $\SMQ_{\Pi_1,\Pi_2}-\SMQ_{\Pi_1,\Pi'_2}$ for
${\rm sp}(\Pi_2)={\rm sp}(\Pi'_2)$ is the ideal $\tilde J$ of
$\WQSym_n\otimes\WQSym_n^*\rightarrow \WQSym_n\otimes\WSym_n^*$
generated by the polynomials
$\W_{\Pi_1}\otimes ({\bf F}^{\Pi_2}-{\bf F}^{\Pi'_2})$.
Since
\begin{equation}
(\WQSym_n\otimes \WQSym_n^*)/_{\tilde J}=\WQSym_n\otimes
\WQSym_n^*/_J=\WQSym_n\otimes \WSym_n^*,
\end{equation}
the result follows.

\medskip
Using the representation of basis elements as pairs of set compositions, one
obtains for two set compositions $\Pi_2$ and $\Pi'_2$ satisfying
${\rm sp}(\Pi_2)={\rm sp}(\Pi'_2)$:
\begin{equation}
\begin{array}{ll}
\Delta(\SMQ_{\Pi_1,\Pi_2}-\SMQ_{\Pi_1,\Pi_2'})=\displaystyle
\sum_{\Pi_1=[\Pi'_1,\Pi''_1]}&
\left(\SMQ_{\std(\Pi'_1),\std(\Pi_2|_{\Pi'_1})}\otimes
\SMQ_{\std(\Pi''_1),\std(\Pi_2|_{\Pi''_1})}\right.\\
&\displaystyle\left.-\SMQ_{\std(\Pi'_1),\std(\Pi'_2|_{\Pi'_1})}\otimes
\SMQ_{\std(\Pi''_1),\std(\Pi'_2|_{\Pi''_1})}\right).
\end{array}
\end{equation}

Since
${\rm sp}(\std(\Pi_2|_{\Pi''_1}))={\rm sp}(\std(\Pi'_2|_{\Pi''_1}))$,
we have
\begin{equation}
(\alpha\otimes\alpha)\circ
\Delta(\SMQ_{\Pi_1,\Pi_2}-\SMQ_{\Pi_1,\Pi_2'})=0.
\end{equation}

Hence, one defines the coproduct $\Delta$ in $\SMCSym$ by making the
following diagram commute
\begin{equation}
\begin{array}{ccc}
\SMQSym&\displaystyle\mathop{\rightarrow}^{\Delta}&
\oplus_n \SMQSym_n\otimes \WQSym_n^*\\
\downarrow \alpha&&\alpha\otimes \alpha'\downarrow\\
\SMCSym&\displaystyle\mathop{\rightarrow}^{\Delta}&
\oplus_n \SMCSym_n\otimes \SMCSym_n^* .
\end{array}
\end{equation}
More precisely, one has
\begin{equation*}
\Delta \SMC_{\Pi_1,\pi_2} = \displaystyle \sum_{\Pi_1=[\Pi'_1,\Pi''_1]}
\SMC_{\std(\Pi'_1),\std(\pi_2|_{\Pi'_1})}\otimes
\SMC_{\std(\Pi''_1),\std(\pi_2|_{\Pi''_1})},
\end{equation*}
where $\pi_2|_{\Pi_1}$ is the set of sets $\{(\pi_2)_i\cap
\bigcup_j(\Pi_1)_j\}_i$ from which one erases the empty sets.

Since $\SMQSym$ is a Hopf algebra, we immediately deduce that
$\Delta$ is an algebra homomorphism from $\SMCSym$ to $\SMCSym\otimes\SMCSym$.
\end{proof}

Note that $\SMRSym$ and $\SMCSym$ have the same Hilbert series, given by
the product of ordered Bell numbers by unordered Bell numbers. This gives one
new example of two different Hopf structures on the same combinatorial set
since $\SMCSym$ is neither commutative nor cocommutative.

Note that the realization of $\SMCSym$ is obtained from the
realization of $\SMQSym$ by quotienting bi-words by the ideal $\mathcal J$
generated by
\begin{equation}
\left\langle \gf{u}{v} \right\rangle
- \left\langle \gf{u}{w} \right\rangle
\end{equation}
where $w$ is obtained from $v$ by permuting its values.

\subsection{Set matrix symmetric functions}

The algebra $\SMSym$ of \emph{set matrix symmetric functions} is the
subalgebra of $\SMCSym$ generated by the polynomials
\begin{equation}
{\rm \SM}_{\pi_1,\pi_2}=\sum_{{\rm sp}(\Pi_1)=\pi_1}\SMC_{\Pi_1,\pi_2}.
\end{equation}

For example,
\begin{equation}
{\rm \SM}_{\{\{1,4\},\{2\},\{3\}\},\{\{1,3,4\},\{2\}\}}
=\sum_{\gf{j_1,\,j_2,\,j_3\text{ distinct}}{k_1 < k_2}}
\left\langle \gf{j_1\,j_2\,j_3\,j_1}{k_1\,k_2\,k_1\,k_1}
 \right\rangle_{\mathcal J}\,.
\end{equation}

\begin{theorem}
\label{TSMSym}
~

\begin{itemize}
\item $\SMSym$ is a co-commutative Hopf subalgebra of $\SMCSym$.
\item $\SMSym$ is isomorphic as an algebra to
${\rm End}_{gr}\WSym=\bigoplus_n \WSym_n\otimes \WSym^*_n$.
\item $\SMSym$ is isomorphic to the quotient of $\SMRSym$ by the ideal
generated by the polynomials $\SMR_{\pi_1,\Pi_2}-\SMR_{\pi_1,\Pi'_2}$
with ${\rm sp}(\Pi_2)={\rm sp}(\Pi'_2)$.
\end{itemize}
\end{theorem}

\begin{proof}
The space $\SMSym$ is stable for the product in $\SMCSym$, as it can be
checked from the realization.
Furthermore, one has
\begin{equation}
\label{DeltaSM}
\begin{array}{rcl}
\Delta \SM_{\pi_1,\pi_2}&=&\sum_{sp(\Pi_1)=\pi_1}\Delta\SMC_{\Pi_1,\pi_2}\\
&=& \sum_{sp(\Pi_1)=\pi_1}\sum_{\Pi_1=[\Pi'_1,\Pi''_1]}
\SMC_{\std(\Pi'_1),\std(\pi_2|_{\Pi'_1})}\otimes
\SMC_{\std(\Pi''_1),\std(\pi_2|_{\Pi''_1})}\\
&=& \sum_{\pi_1=\{\pi'_1,\pi''_1\}}
\SM_{\std(\pi'_1),\std(\pi_2|_{\pi'_1})}\otimes
\SM_{\std(\pi''_1),\std(\pi_2|_{\pi''_1})},
\end{array}
\end{equation}
where $\pi_2|_{\pi_1}$ is the set of sets $\{(\pi_2)_i\cap
\bigcup_j(\pi_1)_j\}_i$ from which one erases the empty sets. The
co-commutativity of $\Delta$ is obvious from~(\ref{DeltaSM}), thus proving the
first part of the theorem.

\medskip
The second part of the theorem is equivalent to the following fact:
$\phi^C(\SM_{\pi_1,\pi_2})= \W_{\pi_1}\otimes F^{\pi_2}.$

\medskip
The proof of the third part is the same as in Proposition~\ref{PSMCSym}.
\end{proof}

\section{Hopf algebras of packed integer matrices}

\subsection{Matrix quasi-symmetric functions}

Let $SA_n$ be the set of set packed matrices such that if one reads
the entries by columns from top to bottom and from left to right,
then one obtains the numbers $1$ to $n$ in the usual order
(see Figure~\ref{SA}).

\begin{figure}[ht]
\begin{equation}
{\left(\begin{array}{cc}
  \st\{1,2\}&\st\es\\\st\es&\st\{4\}\\\st\{3\}&\st\es\end{array}\right)}
\qquad\qquad {\left(\begin{array}{cc}
  \st\{1,3\}&\st\es\\\st\es&\st\{4\}\\\st\{2\}&\st\es\end{array}\right)}
\end{equation}
\caption{\label{SA}An element of $SA$ and an element not in $SA$.}
\end{figure}

Denote by $SA$ the set $SA:=\bigcup_{n\geq 0} SA_n$. One easily sees
that $SA_n$ is in bijection with the packed integer matrices. Indeed,
the bijection $\gimel$ consists in substituting each set of a matrix
by its cardinality. The reverse bijection exists since each integer is the
cardinality of a set, fixed by the reading order of the matrix.
For example,
\begin{equation}
\gimel
\mccc{\{1,2\}&\k8 \emptyset&\k8 \{6\}\\\emptyset&\k8 \{3,4,5\}
&\k8 \{7,8,9,10\}}
=\mccc{2&\k8 0&\k8 1\\0 &\k8 3&\k8 4}.
\end{equation}

Let us consider the subspace $\MQSym'$ of $\SMQSym$ spanned by the elements of
$SA$. For example,
\begin{equation}
\MQ_{\mccc{2&\k8 0&\k8 1\\0 &\k8 3&\k8 4}}:=\SMQ_{\mccc{\{1,2\}&\k8
\emptyset&\k8 \{6\}\\\emptyset&\k8 \{3,4,5\}&\k8 \{7,8,9,10\}}} =
\sum_{\gf{j_1<j_2}{k_1<k_2<k_3}}
 \left\langle
\gf{j_1\,j_1\,j_2\,j_2\,j_2\,j_1\,j_2\,j_2\,j_2\,j_2}%
   {k_1\,k_1\,k_2\,k_2\,k_2\,k_2\,k_3\,k_3\,k_3\,k_3}
\right\rangle\,.
\end{equation}
By definition of the reading order, the product of two elements of $SA$ is a
linear combination of elements of $SA$ and the coproduct of an element of
$SA$ is a linear combination of tensor products of two elements of
$SA$. So

\begin{theorem}
$\MQSym'$ is a Hopf subalgebra of $\SMQSym$ and
it is isomorphic as a Hopf algebra to $\MQSym$.
\end{theorem}

\begin{proof}
As $\MQSym'$ is generated by a set indexed by packed integer matrices, it
is sufficient to check that the product and the coproduct have the same
decompositions than in $\MQSym$. This can be obtained by a straightforward
computation.
\end{proof}

Note that this last theorem gives a realization of bi-words different from the
realization given in~\cite{Hiv}.

\subsection{Matrix half-symmetric functions}

We reproduce the same construction as for set packed matrices. We
define three algebras $\MRSym$ (resp. $\MCSym$, $\MSym$) of packed matrices up
to permutation of rows (resp. of columns, resp. of rows and columns).

\subsubsection{The Hopf algebra $\MRSym$}

Let $\MRSym$ be the subalgebra of $\MQSym$ generated by the polynomials
\begin{equation}
\MR_A:=\sum_B\MQ_{B}
\end{equation}
where $B$ is obtained from $A$ by any permutation of its rows. As for
$\SMRSym$, the realization of $\MRSym$ on bi-words is automatic.
For example,

\begin{equation}
\MR_{\mccc{2&\k8 1&\k8 0\\0&\k8 3&\k8 4}}= \MQ_{\mccc{2&\k8 1&\k8
0\\0&\k8 3&\k8 4}} + \MQ_{\mccc{0&\k8 3&\k8 4\\2&\k8 1&\k8 0}}=
\sum_{\gf{j_1\not=j_2}{k_1<k_2<k_3}}
 \left\langle
\gf{j_1\,j_1\,j_2\,j_2\,j_2\,j_1\,j_2\,j_2\,j_2\,j_2}%
   {k_1\,k_1\,k_2\,k_2\,k_2\,k_2\,k_3\,k_3\,k_3\,k_3}
\right\rangle\,.
\end{equation}

\begin{theorem}
~

\begin{enumerate}
\item $\MRSym$ is a co-commutative Hopf subalgebra of $\MQSym$,
\item $\MRSym$ is also the subalgebra of $\SMRSym$ generated by the
elements $\SMR_A$ where $A$ is any matrix such that each element of the set
composition of its columns is an interval of $[1,n]$.
\end{enumerate}
\end{theorem}

\begin{proof}
From its definition, $\MRSym$ is stable for the product and $\Delta$
maps $\MRSym$ to $\MRSym\otimes\MRSym$. It follows that $\MRSym$ is a Hopf
subalgebra of $\MQSym$. One easily checks the co-commutativity of the
restriction of $\Delta$ to $\MRSym$.

The second part of the theorem amounts to observing that
$\MR_A=\SMR_{\gimel^{-1}A}$.
\end{proof}

\subsubsection{The Hopf algebra $\MCSym$}

We construct the algebra $\MCSym$ as the quotient of $\MQSym$ by the
ideal generated by the polynomials ${\rm \MQ}_A-{\rm \MQ}_B$ where
$B$ can be obtained from $A$ by a permutation of its columns.

\begin{theorem}
~

\begin{enumerate}
\item $\MCSym$ is a commutative Hopf algebra,
\item $\MCSym$ is isomorphic as a Hopf algebra to the subalgebra of $\SMCSym$
generated by the elements $\SMC_A$, such that each element of the set
partition of its columns is an interval of $[1,n]$.
\end{enumerate}
\end{theorem}

\begin{proof}
The proof of the first part of the theorem is almost the same as
Proposition~\ref{PSMCSym}(1).

The dimensions are the same, so it is sufficient to check that the
product and the coproduct have the same decomposition in both algebras.
\end{proof}

\subsubsection{Dimensions of $\MRSym$ and $\MCSym$}

The dimension of the homogeneous component of degree $n$ of $\MRSym$
or $\MCSym$ is equal to the number of packed matrices with sum of entries
equal to $n$, up to a permutation of their rows.

Let us denote by $\PMuR(p,q,n)$ the number of such $p\times q$
matrices. One has obviously
\begin{equation}\label{dimMR1}
\dim \MRSym_n=\sum_{1\leq p,q\leq n}\PMuR(p,q,n).
\end{equation}
The integers   $\PMuR(p,q,n)$ can be computed through the induction
\begin{equation}
 \PMuR(p,q,n)=\MuR(p,q,n)-\sum_{1\leq k,l\leq
 p,q}\left(\gf{q}{l}\right)\PMuR(k,l,n),
\end{equation}
where $\MuR(p,q,n)$ is the number of $p\times q$ possibly unpacked matrices
with sum of entries equal to $n$, up to a permutation of their rows.

Solving this induction and substituting it in equation (\ref{dimMR1}), one
gets
 \begin{equation}
\dim \MRSym_n=\sum_{i=1}^{n+1}(-1)^{n-i}T_{n+1,i+1}\MuR(n,i,n),
 \end{equation}
where
\begin{equation}
T_{n,k}=\sum_{j=0}^{n-k}(-1)^{n-k-j}\left(\gf{j+k-1}{j}\right)
\end{equation}
that is the number of minimum covers of an unlabeled $n$-set that cover $k$
points of that set uniquely (sequence A056885 of \cite{Slo}).
The generating series of the $T_{n,k}$ is
\begin{equation}
\sum_{i,j}T_{i,j}x^iy^j={\frac {1-x}{\left (1+x\right )\left (1-x-xy\right )}}.
\end{equation}
The integer $\MuR(n,i,n)$, computed via the  P\'olya enumeration
theorem, is the coefficient of $x^n$ in the cycle index
$Z(G_{n,i})$, evaluated over the alphabet $1+x+\cdots+x^n+\cdots$,
of the subgroup $G_{n_,i}$ of $\S_{in}$ generated by the
permutations $\sigma. \sigma[n]. \sigma[2n].\ \cdots\ . \sigma[in]$
for $\sigma\in \S_n$ (here $.$ denotes the concatenation).

This coefficient is also the number of partitions $N_{n,i}$ of $n$
objects with $i$ colors whose generating series is
\begin{equation}
\sum_{n}N_{n,i}x^n
=\prod_k
 \left(\frac{1}{1-x^k}\right)^{\left(\gf{i+k}{i}\right)}.
\end{equation}
Hence
\begin{proposition}
 \begin{equation}
\dim \MRSym_n=\sum_{i=1}^{n+1}(-1)^{n-i}T_{n+1,i+1}N_{n,i}.
 \end{equation}
\end{proposition}

The first values are
\begin{equation}
\begin{split}
{\rm Hilb}(\MRSym)={\rm
Hilb}(\MCSym)=&
\; 1+t+4\,t^2+16\,t^3+76\,t^4+400\,t^5+2356\,t^6+15200\,t^7\\
&+106644\,t^8 + 806320\,t^9+6526580\,t^{10}+\cdots
\end{split}
\end{equation}

\subsection{Matrix symmetric functions}

The algebra $\MSym$ of \emph{matrix symmetric functions} is the
subalgebra of $\MCSym$ generated by

\begin{equation}
 \M_A:=\sum_{B}\MR_{B}, 
\end{equation}
where $B$ is obtained from $A$ by any permutation of its rows.

\begin{theorem}
~

\begin{enumerate}
\item $\MSym$ is a commutative and co-commutative Hopf subalgebra of
$\MCSym$.
\item $\MSym$ is isomorphic as a Hopf algebra to the subalgebra of $\SMSym$
generated by the elements $\SM_A$, such that each element of the set
partition of its columns is an interval of $[1,n]$.
\item $\MSym$ is isomorphic to the quotient of $\MRSym$ by the ideal generated
by the polynomials ${\rm \MR}_A-{\rm \MR}_B$ where $B$ can be obtained from
$A$ by a permutation of its columns.
\end{enumerate}
\end{theorem}

\begin{proof}
The proof follows the same lines as the proof of Theorem~\ref{TSMSym}.
\end{proof}

From all the previous results, we deduce that the following diagram commutes

\begin{equation*}
\newdir{ >}{{}*!/-4mm/@{>}}
\xymatrix@R=10mm@L=10mm{& \SMQSym \ar@{<-<}[dl] \ar@{->>}[ddr]\ar@{<-<}[rrr]&&&\MQSym\ar@{<-<}[dl] \ar@{->>}[ddr]\\
\SMRSym\ar@{->>}[ddr]\ar@{<-<}[rrr]& &&\MRSym\ar@{->>}[ddr]\\
&&\SMCSym\ar@{<-<}[dl]\ar@{<-<}[rrr]&&&\MCSym\ar@{<-<}[dl]\\
  &\SMSym\ar@{<-<}[rrr]&&&\MSym}
\end{equation*}

\section{Dendriform structures over $\SMQSym$}

\subsection{Tridendriform structure}

A {\em tridendriform algebra} is an associative algebra whose multiplication
can be split into three operations
\begin{equation}
x\cdot y = x\prec y + x\circ y + x\succ y\,,
\end{equation}
where $\circ$ is associative, and such that
\begin{equation}
(x\prec y)\prec z = x\prec (y\cdot z)\,,\ \
(x\succ y)\prec z = x\succ (y\prec z)\,,\ \
(x\cdot y)\succ z = x\succ (y\succ z)\,,\ \
\end{equation}
\begin{equation}
(x\succ y)\circ z = x\succ (y\circ z)\,,\ \ \
(x\prec y)\circ z = x\circ (y\succ z)\,,\ \ \
(x\circ y)\prec z = x\circ (y\prec z)\,.
\end{equation}

\subsection{Tridendriform structure on bi-words}

One defines three product rules over bi-words as follows:
\begin{enumerate}
\item
$\left\langle \gf{u_1}{v_1}\right\rangle \prec
 \left\langle \gf{u_2}{v_2}\right\rangle
=\left\langle \gf{u_1}{v_1}\right\rangle \star
 \left\langle \gf{u_2}{v_2}\right\rangle
 \mbox{\ \ if $\max(u_1)>\max(u_2)$,\ and $0$ otherwise.}$
\item
$\left\langle \gf{u_1}{v_1}\right\rangle \,\,\circ\,\,
 \left\langle \gf{u_2}{v_2}\right\rangle
=\left\langle \gf{u_1}{v_1}\right\rangle \star
 \left\langle \gf{u_2}{v_2}\right\rangle
 \mbox{\ \ if $\max(u_1)=\max(u_2)$,\ and $0$ otherwise.}$
\item
$\left\langle \gf{u_1}{v_1}\right\rangle \succ
 \left\langle \gf{u_2}{v_2}\right\rangle
=\left\langle \gf{u_1}{v_1}\right\rangle \star
 \left\langle \gf{u_2}{v_2}\right\rangle
 \mbox{\ \ if $\max(u_1)<\max(u_2)$,\ and $0$ otherwise.}$
\end{enumerate}

\begin{proposition}
The algebra of bi-words endowed with the three product rules $\prec, \circ$,
and $\succ$ is a tridendriform algebra.

Moreover, $\SMQSym$ is stable by those three rules.
More precisely, one has:
\begin{equation}
\label{precSM}
\SMQ_{\left\langle \gf{u}{\bf v}\right\rangle}
\prec\SMQ_{\left\langle \gf{u'}{\bf v'}\right\rangle}
=\displaystyle
 \sum_{\gf{w=x.y\in u\star_W u'}{|x|=|u|;\ \max(y)<\max(x)}}
 \SMQ_{\left\langle \gf{w}{{\bf v}{\bf v'}[\max(v)]}\right\rangle},
\end{equation}
\begin{equation}
\label{circSM}
\SMQ_{\left\langle \gf{u}{\bf v}\right\rangle}
\circ\SMQ_{\left\langle \gf{u'}{\bf v'}\right\rangle}
=\displaystyle
 \sum_{\gf{w=x.y\in u\star_W u'}{|x|=|u|;\ \max(y)=\max(x)}}
 \SMQ_{\left\langle \gf{w}{{\bf v}{\bf v'}[\max(v)]}\right\rangle},
\end{equation}
\begin{equation}
\label{succSM}
\SMQ_{\left\langle \gf{u}{\bf v}\right\rangle}
\succ\SMQ_{\left\langle \gf{u'}{\bf v'}\right\rangle}
=\displaystyle
 \sum_{\gf{w=x.y\in u\star_W u'}{|x|=|u|;\ \max(y)>\max(x)}}
 \SMQ_{\left\langle \gf{w}{{\bf v}{\bf v'}[\max(v)]}\right\rangle}.
\end{equation}

So $\SMQSym$ is a tridendriform algebra.
\end{proposition}

\begin{proof}
The first part of the proposition amounts to checking the compatibility
relations between the three rules. It is immediate.

The stability of $\SMQSym$ by any of the three rules and the product relations
are also immediate: the bottom row can be any word whose packed word is
${\bf v}{\bf v'}[\max(v)]$, so that the formula reduces to a formula on the
top row which is equivalent to the same computation in $\WQSym$
(see~\cite{NT}). The compatibility relations automatically follow from their
compatibility at the level of bi-words.
\end{proof}

\begin{corollary}
$\SMCSym$, $\MQSym$ and $\MCSym$ are tridendriform.
\end{corollary}

\begin{proof}
As in the case of $\SMQSym$, one only has to check that the algebras are
stable by the three product rules since the compatibility relations
automatically follow.

The case of $\SMCSym$ is direct since Formulas~(\ref{precSM}),~(\ref{circSM}),
and~(\ref{succSM}) have only the word ${\bf v}{\bf v'}[\max(v)]$ in the
bottom row of their basis elements. For the same reason, the case of $\MCSym$
directly follows from the case of $\MQSym$.
The case of $\MQSym$ is the same as the case of $\SMQSym$ itself.
\end{proof}

\subsection{Bidendriform structures}

Let us define two product rules $\ll=\prec$ and $\gg=\circ+\succ$ on bi-words.
We now split the non-trivial parts of the coproduct of the $\SMQ$ of
$\SMQSym$, as
\begin{equation}
\Delta_\ll(\SMQ_{A})=
\displaystyle
\sum_{\gf{A=\left(\gf{B}{C}\right),\ A\neq B,C}{\max(B)=\max(A)}}
 \SMQ_{\std(B)}\otimes \SMQ_{\std(C)}.
\end{equation}
\begin{equation}
\Delta_\gg(\SMQ_{A})=
\displaystyle
\sum_{\gf{A=\left(\gf{B}{C}\right),\ A\neq B,C}{\max(C)=\max(A)}}
 \SMQ_{\std(B)}\otimes \SMQ_{\std(C)}.
\end{equation}

Let us recall that under certain compatibility relations between the two parts
of the coproduct and other compatibility relations between the two product
rules $\ll:=\prec$ and $\gg:=\circ+\succ$ defined by Foissy~\cite{Foissy},
we get bidendriform bialgebras.

\begin{theorem}
$\SMQSym$ is a bidendriform bialgebra.
\end{theorem}

\begin{proof}
The co-dendriform relations, the one concerning the two parts of the
coproduct, are easy to check since they only amount to knowing which part of a
matrix cut in three contains its maximum letter.

The bi-dendriform relations are more complicated but reduce to a careful check
that any part of the coproduct applied to any part of the product only brings
a limited amount of disjoint cases. Let us for example check the relation
\begin{equation}
\label{bidend2}
\Delta_\droitdend  (a\gaudend b) =
      a'b'_\droitdend \!\otimes\! a''\!\gaudend \!b''_\droitdend 
\,+\, a'\!\otimes\! a''\!\gaudend \!b 
\,+\, b'_\droitdend \!\otimes\! a\!\gaudend \!b''_\droitdend \,,
\end{equation}
where the pairs $(x',x'')$ (resp. $(x'_\gaudend ,x''_\gaudend )$ and
$(x'_\droitdend,x''_\droitdend)$) correspond to all possible elements
occurring in $\overline\Delta x$ (resp. $\Delta_\gaudend x$ and
$\Delta_\droitdend  x$), summation signs being understood (Sweedler's
notation). 

First, the last row of all elements in $(a\gaudend b)$ only contain elements
of $a$.
Since by application of $\Delta_\droitdend$, the maximum of $b$ has to go in
the right part of the tensor product, this means that there has to be also
elements coming from $a$ in this part of the tensor product.
Now, the elements of $\Delta_\droitdend(a\gaudend b)$ where all elements of
$b$ are in the right part of the tensor product, are obtained, for the left
part by elements coming from the top rows of $a$ and for the right part by
elements coming from the other rows of $a$ multiplied by $b$ in such a way
that the last row only contains elements of $a$, hence justifying
the middle term $a'\!\otimes\! a''\!\gaudend \!b$.

If both components of $\Delta_\droitdend(a\gaudend b)$ contain elements coming
from $b$, then the left part cannot contain the maximum element of $b$ (hence
justifying the $b'_\droitdend$ and $b''_\droitdend$, the left part being
multiplied by elements coming from $a$ if any (this is the difference between
the first and the third term of the expansion of
$\Delta_\droitdend (a\gaudend b)$), the right part being multiplied by
$\gaudend$ with elements coming from $a$ since the last row must contain
elements coming from $a$.
\end{proof}

Recall that $\SMCSym$ is the quotient of $\SMQSym$ by the ideal generated by
$\SMQ_A - \SMQ_B$ where $A$ and $B$ are the same matrices up to a permutation
of their columns, the row containing the maximum element is the same for any
element of a given class, so that the left coproduct and the right coproduct
are compatible with the quotient.
Moreover, the left and right coproduct are internal within $\MQSym$, so that

\begin{corollary}
$\SMCSym$ and $\MQSym$ are bidendriform sub-bialgebras of $\SMQSym$.
\end{corollary}

\begin{proof}
We already know that $\SMCSym$ and $\MQSym$ are dendriform subalgebras of
$\SMQSym$ since they are tridendriform subalgebras of this algebra. The
compatibility relations come from the compatibility relations on $\SMQSym$, so
that there only remains to check that the coproduct goes from $X$ to
$X\otimes X$, where $X$ is either $\SMCSym$ or $\MQSym$. This is an easy
computation.
\end{proof}

\begin{corollary}
$\SMQSym$, $\SMCSym$, and $\MQSym$ are free, cofree, self-dual Hopf algebras
and their primitive Lie algebras are free.
\end{corollary}

\begin{proof}
This follows from the characterization of bidendriform bialgebras done by
Foissy~\cite{Foissy}.
\end{proof}

\subsection{A realization of $LDiag(q_c,q_s)$ on bi-words}

Let us define a two-parameter generalization of the algebra $\MQSym$. For this
purpose, consider bi-words with parameter-commuting bi-letters depending on
the bi-letters as follows:

\begin{equation}
\begin{split}
\left\langle \gf{yx}{zt} \right\rangle
=q_c \left\langle \gf{xy}{tz} \right\rangle
&\mbox{ if }y>x, \\
\left\langle \gf{xx}{zt} \right\rangle
=q_s \left\langle \gf{xx}{tz} \right\rangle
&\mbox{ if }z<t.\\
\end{split}
\end{equation}

Let us now define the realization as a sum of bi-words of a packed integer
matrix with $p$ rows and $q$ columns:
\begin{equation}
{\bf LD}_M :=
     \sum_{\gf{j_1<\dots<j_p}{k_1<\dots<k_q}}
        \prod_{a=1}^p \prod_{b=q}^1
     \left\langle \gf{j_a}{k_b} \right\rangle^{M_{ab}}.
\end{equation}

For example,
\begin{equation}
{\bf LD}_{\mcc{3&\k8 5\\1&\k8  3}}
= \sum_{\gf{j_1<j_2}{k_1<k_2}}
   \left\langle
    \gf{j_1^5\,j_1^3\,j_2^3\,j_2}{k_2^5\,k_1^3\,k_2^3\,k_1}
   \right\rangle.
\end{equation}

We then have
\begin{theorem}
~

\begin{itemize}
\item The subspace spanned by the $\bf LD$ has a structure of associative
algebra.
Moreover, the matrices indexing the product ${\bf LD}_A {\bf LD}_B$ are
equal to the matrices appearing in $\M_A \M_B$ in $\MQSym$, and the
coefficient of ${\bf LD}_C$ in this product is a monomial
$q_s^{x(A,B,C)}q_c^{y(A,B,C)}$ computed as follows:
let us call \emph{left} the part of $C$ coming from $A$ and \emph{right} the
part of $C$ coming from $B$. Then 
\begin{equation}
x(A,B,C) = \sum_{r \text{ row of $C$}}
\left(\sum_{i\in\text{ left(r)}}i\right)
\left( \sum_{j\in\text{ right(r)}}j\right).
\end{equation}
\begin{equation}
y(A,B,C) = \sum_{r<r' \text{ rows of $C$}}
\left(\sum_{i\in\text{ left(r)}}i\right)
\left(\sum_{j\in\text{ right(r')}}j\right).
\end{equation}
\item The specialization $q_s=q_c=1$ gives back $\MQSym$.
\end{itemize}
\end{theorem}

\begin{proof}
Since each partially commuting bi-word has only one expression such that the
top row is weakly increasing and the bottom row is weakly decreasing at the
spots where the top row is constant, we can define without ambiguity the
\emph{canonical} element of a bi-word. The set of canonical elements is in
bijection with integer matrices.

Now, if two bi-words appearing in a product of two ${\bf LD}$ have canonical
elements whose corresponding matrices have the same packed matrix,
they follow exactly the same rewriting steps to get to their canonical
element. So in particular, the product of two ${\bf LD}$ decomposes as a
linear combination of ${\bf LD}$.
Moreover, since the product on bi-words is associative and compatible with the
partial commutations, then so is the product of the ${\bf LD}$, hence proving
that they span an algebra.

By definition of the realization of $\MQSym$ on bi-words, $\MQSym$ is
obtained from this algebra by specifying $q_c=q_s=1$, that is, replacing
partially parameter-commuting bi-letters by partially commuting bi-letters, so
that the matrices appearing in a product of two ${\bf LD}$ are the same as the
matrices appearing in the product of the same packed matrices in $\MQSym$.
Finally, the coefficient of a given matrix is obviously a monomial in $q_s$
and $q_c$ and the powers of $q_s$ and $q_c$ are straightforward from the
definition of the commutations: a bi-letter of the right has to exchange with
any bi-letter of the left whose top value is greater than or equal to its top
value. Each exchange amounts either to multiplying by $q_c$ if those values
differ, or to multiplying by $q_s$ is they are equal. This is equivalent to
the formulas of the statement.
\end{proof}

For example, one has:
\begin{equation}
\begin{split}
{\bf LD}_{\mccc{ 2 &\k8 0 \\ 1 &\k8 4}}
\star {\bf LD}_{\left(1\right)}
=& {\bf LD}_{\mccc{
        2 &\k8 0 &\k8 0 \\
        1 &\k8 4 &\k8 0 \\
        0 &\k8 0 &\k8 1 }}
+ q_s^5\, {\bf LD}_{\mccc{
                  2 &\k8 0 &\k8 0 \\
                  1 &\k8 4 &\k8 1 }}
+ q_c^5\, {\bf LD}_{\mccc{
                  2 &\k8 0 &\k8 0 \\
                  0 &\k8 0 &\k8 1 \\
                  1 &\k8 4 &\k8 0}}
\\
&+ q_c^5q_s^2\, {\bf LD}_{\mccc{
                      2 &\k8 0 &\k8 1 \\
                      1 &\k8 4 &\k8 0}}
+ q_c^7\, {\bf LD}_{\mccc{
                  0 &\k8 0 &\k8 1 \\
                  2 &\k8 0 &\k8 0 \\
                  1 &\k8 4 &\k8 0 }}
\end{split}
\end{equation}

since
\begin{equation}
\begin{split}
\sum_{\gf{j_1<j_2}{k_1<k_2}}
\left\langle
\gf{j_1\,j_1\,j_2\,j_2\,j_2\,j_2\,j_2}{k_1k_1k_2k_2k_2k_2k_1}
\right\rangle
\star \sum \left\langle \gf{j}{k} \right\rangle
= &
\sum_{\gf{j_1<j_2<j_3}{k_1<k_2<k_3}}
   \left\langle \gf{j_1^2\,j_2^4\, j_2\, j_3}{k_1^2\,k_2^4\,k_1\,k_3}
   \right\rangle
+
\sum_{\gf{j_1<j_2(=j_3)}{k_1<k_2<k_3}}
   q_s^5\left\langle \gf{j_1^2\,j_2\,j_2^4j_2}{k_1^2\,k_3\,k_2^4\,k_1}
  \right\rangle\\
&+
\sum_{\gf{j_1<j_3<j_2}{k_1<k_2<k_3}}
   q_c^5\left\langle \gf{j_1^2\,j_3\,j_2^4\,j_2}{k_1^2\,k_3\,k_2^4\,k_1}
   \right\rangle
+
\sum_{\gf{j_1(=j_3)<j_2}{k_1<k_2<k_3}}
  q_s^2q_c^5\left\langle \gf{j_1\,j_1^2\,j_2^4\,j_2}{k_3\,k_1^2\,k_2^4\,k_1}
            \right\rangle\\
&+
\sum_{\gf{j_3<j_1<j_2}{k_1<k_2<k_3}}
 q_c^7
 \left\langle \gf{j_3\,j_1^2\,j_2^4\,j_2}{k_3\,k_1^2\,k_2^4\,k_1}
 \right\rangle.
\end{split}
\end{equation}


\end{document}